\documentclass[11pt]{amsart}\usepackage{amsfonts} \usepackage{amscd}
\usepackage{amssymb}
\usepackage{amscd}

\newcommand{\be}{\begin{equation}}
\newcommand{\ee}{\end{equation}}

\newcommand{\ring}[1]{\mathbb{#1}}
 
\newcommand{\Z}{\ring{Z}} \newcommand{\R}{\ring{R}}

\newcommand{\N}{\ring{N}}

\newcommand{\ldva}{\ell^2}
\newcommand{\co}{{\mathcal C}}
\newcommand{\ba}{{\mathcal B}}
\newcommand{\lt}{\lambda B}

\title[Common hypercyclic vectors]{Common hypercyclic vectors for multiples of
backward shift}
\author{Evgeny Abakumov \& Julia Gordon}

\address{Equipe d'Analyse et de Math\'ematiques Appliqu\'ees,
Universit\'e de Marne-la-Vall\'ee, 5 Boulevard Descartes,
Champs-sur-Marne, 77454 Marne-la-Vall\'ee CEDEX 2, France}
\email{abakumov@math.univ-mlv.fr}
\address{Department of Mathematics, University of Michigan, Ann Arbor
  MI 48109-1109, USA}
\email{julygord@umich.edu
}

\begin{document}
\theoremstyle{plain}
\newtheorem*{thm}{Theorem}
\newtheorem*{lem}{Lemma}

\theoremstyle{definition}
\newtheorem*{rem}{Remark}

\footnotetext[1]{Keywords: hypercyclic vector, backward shift}
\footnotetext[2]
{2000 {\it Mathematics Subject Classification}: 47A16, 47B37}

\begin{abstract}

We prove that  $\ell^2$ contains  a dense set of vectors which are
hypercyclic simultaneously for all multiples of the
backward shift operator by
constants of absolute value greater than $1$.

\end{abstract}
\maketitle

\setcounter{section}{-1}

\section{Introduction}

The focus
of our study in this paper is the
hypercyclic behavior (density of an orbit) for scalar multiples of 
the backward shift operator.

For a linear topological space $X$, a continuous  linear operator 
$T\colon X\to X$
is called {\it hypercyclic} if there exists a vector $x\in X$ such that
its orbit   $\{T^n x,  n \geq 0\}$ is dense in $X$ (such a vector
$x$ is called a {\it hypercyclic vector} for the operator $T$).

The study of the phenomenon of hypercyclicity originates in the papers
by Birkhoff, 1929 \cite{B} and
MacLane, 1952  \cite{M} that show, respectively, that the operators of
translation and differentiation, acting on the space of entire functions,
are hypercyclic.
The first example of a hypercyclic operator acting on a
Banach space was given by  Rolewicz in
1969 \cite{R}. The example is the backward shift operator on $\ell^p$, 
scaled by a  constant greater than $1$.

Later on, 
starting from the eighties, 
the subject of hypercyclicity has been widely 
explored.
We refer the reader  to the paper by Godefroy and Shapiro \cite{GS}
and the recent survey by Grosse-Erdmann \cite{G}
for information  and references.

We will concentrate on the space $\ldva$ of all complex valued
square summable sequences.
The backward shift operator $B$ on this space is defined by
$$ B(a_0, a_1, a_2, \dots) = (a_1, a_2, a_3, \dots).
$$
Since the operator $B$ is a contraction, $B$ 
itself cannot be hypercyclic. The above mentioned
 result due to Ro\-le\-wicz 
states that the operator $zB$ has a hypercyclic vector
whenever $|z|>1$.    
The question that we answer in this paper is
whether  it is possible to find a {\it single} vector
$f\in\ldva$ such that its orbit
 is dense
in $\ldva$ for {\it all} operators  $z B$ with $|z|>1$. The answer is {\it yes}:

\begin{thm}\label{th}
Let $B$ be the  backward shift acting on $\ell^2$.
There exists a common hypercyclic vector for the operators
$zB, \, |z|>1$.
\end{thm}
 For a continuous
linear operator  $T:X\to X$ acting on a Banach  space $X$, the notation
 $HC(T)$ stands for the set of all vectors in $X$ that are hypercyclic for
$T$.
 It turns out that the set $HC(T)$ is residual (i.e. dense 
 $G_\delta$)
whenever it is  nonempty,
see, e.g., \cite{GS}.
This fact together with Baire Category Theorem immediately implies that
if $\{T_\lambda, \lambda \in \Lambda\}$ is a family of hypercyclic operators,
then $\cap_{\lambda \in \Lambda} HC(T_\lambda)$ is residual (in particular,
nonempty) whenever $\Lambda$ is {\it countable}.

For {\it uncountable} families of operators the
standard Baire categories arguments fail, and this is what
makes our question interesting. At the moment we do not know of any other
non-trivial examples
of uncountable families of operators with common hypercyclic vectors.
The question  about  the set
$\cap_{|z|>1}HC(zB)$ was asked  by
Salas \cite{S}. Our theorem
states that this set is nonempty, and hence dense 
(it suffices to notice that if $x$ is a hypercyclic vector
for an operator $T$, then each element of its orbit $\{T^nx,\,n\ge0\}$ is also 
hypercyclic for $T$).  We do not know
if  the set of common hypercyclic vectors is large or small (say, in the 
sense of Baire categories).

Let us mention here that analogous problems for {\it cyclic} vectors
(those whose orbit {\it spans} a dense subset in the space) have been
treated for some natural families of operators. In particular, there exists
a common cyclic vector for
adjoint multiplications on some spaces of analytic functions (see \cite{W},
\cite{C1},
\cite{BS}), and for the collection of linear partial
differential operators of finite order  with constant coefficients
 (see \cite{C2}).

We will use the following notations:
$\N=\{1,2,\ldots\}$ -- the set of natural numbers;
$\R_+$ -- the set of positive real numbers;
$[r]$ -- the integral part of a real number, $[r]=\max\{n\in\Z\mid n\leq r\}$;
$B$ -- the backward shift operator on the space $\ldva$, as defined above;
by $\|\cdot\|$ we will always mean $\ldva$-norm.
%

{\bf Acknowledgement.}
We are grateful to
 Alexander Borichev for reading the manuscript and for his
helpful suggestions.

\section{Algorithm}
We start with a review of the main constituents of our construction.
For the simplicity of exposition,
we construct a common hypercyclic vector for the family of real multiples
$\{\lambda B$, $\lambda>1\}$ of the operator $B$.
Then, in Section \ref{complex}, we describe how to modify the construction
so that it would give a common hypercyclic vector for all complex
multiples $zB$, $|z|>1$.

\subsection{Preliminary remarks}
The first obvious observation is that if $f$ is any vector in $\ldva$ and
$n\in\N$ is fixed,
then all the vectors $(\lambda B)^n f$, $\lambda>0$,  have the same
direction as the vector $B^n f$. In particular, if we want $f$ to be
a common hypercyclic vector for the family $(\lambda B)_{\lambda >1}$, then
the sequence $B^n f$ should ``approximate all directions''. 
To better define this requirement, we introduce a countable
family of cones generated by a set of diminishing in size
 balls in $\ldva$ centered at the elements of
a dense set in $\ldva$. One of the requirements on a common hypercyclic
vector $f$  is
that $B^n f$ falls into every cone sufficiently often.

To specify what ``sufficiently often'' means,
we have to look closely at the norms
of the vectors $B^n f$, where $f$ is a common hypercyclic vector.
Let us fix one of the cones $\co$ from the
previous paragraph, and let $\{n_k\}_{k\geq 1}$ be the set of all integers $n$
such that $B^n f\in \co$.  Then, naturally, for all $\lambda>1$ we have
$(\lambda B)^n f\in\co$.
Let $\ba$ be a ball inside $\co$.
Then for every $\lambda>1$
there exists
$k\in\N$ such that $\lambda^{n_k}B^{n_k}f\in \ba$. In particular, this implies
that  for any $\lambda>1$ there should exist an integer $k$ such that
the  $\ldva$-norm of the vector  $\lambda^{n_k}B^{n_k}f$ belongs to a
certain fixed interval $]a-\delta, a+\delta[\subset\R$.
This turns out to be a very restrictive condition
on the  sequence of the norms $\{\|B^{n_k} f\|\}_{k\geq 1}$,
as illustrated by the following
statement.

\begin{lem}\label{equi}
Let $\{n_k\}_{k\geq 1}$ be a fixed increasing sequence of positive integers,
and let $a\in\R_+$.
Then the following conditions on a sequence $\{r_k\}_{k\geq 1}$ of real
 positive
numbers are equivalent:
\newcounter{cond}
\begin{list}%
{\arabic{cond})}{\usecounter{cond}\setlength{\leftmargin}{0.0in}}
\item
For every   $\lambda>1$, $\epsilon >0$,
there exists $k\in\N$ such that $|\lambda^{n_k}r_k-a|<\epsilon$.
\item
Let $x_k=-\frac{\ln r_k}{n_k}$. Then
for any $\delta>0$, the
$\frac{\delta}{n_k}$-neighbourhoods
of the points $\frac{\ln a}{n_k}+x_k$ cover $\R_+$.
\end{list}
In particular, the sequence $\{\lambda^{n_k}r_k\}_{k\geq 1}$
is dense in $\R_+$ for any
$\lambda>1$  if and only if
 the condition (2) holds for a dense set of $a$'s.
\end{lem}
Here the sequence $\{r_k\}_{k\geq 1}$ plays the role of
 $\{\|B^{n_k} f\|\}_{k\geq 1}$.

{\bf Proof.} The proof is straightforward. \qed

\begin{rem} The lemma shows in particular that
it is not always possible to find such a sequence
$\{r_k\}_{k\geq 1}$ for a given sequence
$\{n_k\}_{k\geq 1}$. For example,
($2$) implies that the Lebesgue measure of the union
of  $\frac{\delta}{n_k}$-neighbourhoods of
certain points is infinite;
therefore, if $\{r_k\}_{k\geq 1}$ exists, then
the series $\sum_{k\ge1}\frac1{n_k}$ necessarily
diverges.
\end{rem}

It follows from this remark 
 that if $f$ is a common hypercyclic vector, and $\{n_k\}_{k\ge 1}$
is the sequence of all numbers such that the iterates $B^{n_k}f$
belong to a given cone, then 
the series $\sum_{k\ge 1}\frac1{n_k}$ must diverge,
for every fixed cone. For this condition to hold,
the iterates of the hypercyclic vector
are required to
``visit'' every cone of our infinite family very often.

\subsection{Main idea}\label{mi}
Here is the algorithm that we follow in order to construct the common
hypercyclic vector:
\newcounter{alg}
\begin{list}%
{\arabic{alg})}{\usecounter{alg}\setlength{\leftmargin}{0.3in}}
\item\label{vs}
 Fix a dense countabe set in $\ldva$.
\item\label{cones}
Construct a  sequence of cones generated by balls  centered at the
vectors of this set with radii going to
$0$.  The sequence of all these cones will be denoted
by $\{\co_l\}_{l=1}^{\infty}$.
\item\label{ji}
Fix a function $j\colon\{n\in \N, n\geq 20\} \to\N$ that takes every
 positive integral value
sufficiently often, and possesses certain additional properties.
\item\label{seq}
Construct a  sequence $\{M_n\}_{n\geq 20}$ of natural numbers simultaneously
with
 a sequence $\{r_n\}_{n\geq 20}$  of positive real numbers.
The sequence $\{M_n\}_{n\ge 20}$ will be the sequence of numbers of
 the iterates of the
prospective common hypercyclic vector, over which we have complete control
(see the next step).
This sequence will fall into a disjoint union of
infinitely many subsequences according to the values of the function $j$:
the $l$-th subsequence will consist of the elements of $\{M_n\}$ with
$j(n)=l$.
If we set $\{n_k\}_{k\ge 1}$ to be any of these subsequences
 without a few initial terms,
 and
$\{r_k\}_{k\ge 1}$ to be the corresponding
subsequence of the sequence $\{r_n\}_{n\ge 20}$,
then $\{n_k\}_{k\ge 1}$ and
$\{r_k\}_{k\ge 1}$
satisfy the condition ($2$) of the Lemma 
and some additional requirements.

\item\label{vect}
Construct a  vector $f\in\ldva$ such that $\|B^{M_k} f\|=r_k$ for
all $k\geq 20$, and such that for any $l\ge1$, $B^{M_k} f$ belongs to the cone
 $\co_{l}$
for  all sufficiently large $k$ verifying $j(k)=l$.
\end{list}

It is roughly clear from the previous subsection that such a vector
will be a common hypercyclic vector for the operators $\lt$, $\lambda>1$.
We give a rigorous proof of this statement in Section $3$.
Section  $2$ is devoted to Step (\ref{seq}) of the algorithm, that is,
to the construction of the sequences $\{M_k\}$ and $\{r_k\}$.

\subsection{Notation} Here we carry out the Steps (\ref{vs}),
(\ref{cones}), and (\ref{ji}) of the algorithm.
\subsubsection{}\label{vectorsandcones}
We call a vector $v \in\ldva$ {\it finite}, if all its coordinates,
except for finitely many, equal $0$.

Let $\{v_l\}_{l\geq 1}$ be a dense subset of
$\ldva$, such that for all $l\in\N$  the vector $v_l$ is nonzero
and finite, and
put $\alpha_l=\ln\|v_l\|$.
Clearly, it is possible to choose the sequence $\{v_l\}_{l\geq 1}$
in such a way that
the real numbers $\alpha_l$ satisfy the following condition:
\begin{equation}\label{alfa}
|\alpha_{l+1}-\alpha_l|\leq 1 {\text {\ for all \ }}l\geq 1;
\quad \alpha_1=\alpha_2=0.
\end{equation}
 We fix a sequence  of positive real numbers
$\{\epsilon_l\}_{l\geq 1}$ such that
$\epsilon_l<\|v_l\|$ and $\epsilon_l\to 0$ as $l\to\infty$, and
define the cone $\co_l$ by $\co_l=\{v\in\ldva\mid \exists \mu>0\colon
\|\mu v-v_l\|<\epsilon_l\}$.

\subsubsection{}\label{dzhi}
We define the function $j\colon \{n \in\N, n\geq20\}\to\N$ in the following
way:  let  $j(k)=k+1-[n\ln\ln n]$, where $n$ is the unique integer satisfying
the inequality
$[n\ln\ln n]\leq k < [(n+1)\ln\ln(n+1)]$. (Notice that $j(20)=1$.)
This function takes every positive integral value infinitely many times:
its domain is subdivided into  intervals, on each of which
$j(k)$ increases with unit step  from the value $1$ at the beginning of the
interval. The constant 20 appears 
in the construction due to 
the inequality $\ln\ln 20 >1$ which will be used below.


We will use the following  property of the function $j(k)$ which is easy to
check:
\begin{equation}\label{jin}
j(k)<[\ln\ln k]+3 {\text{\ for\ }} k\geq 20.
\end{equation}
Now we are ready to begin the actual construction.

\section{Main construction}
Here we carry out Step (\ref{seq}) of the algorithm declared in the
previous section.
Let $\{\alpha_l\}_{l\geq 1}$ be the sequence of real numbers from Section
\ref{vectorsandcones}.

We construct two sequences: $\{M_k\}_{k\geq 20}$ of positive integers
 and \newline
$\{x_k\}_{k\geq 20}$ of real numbers, that possess the following
properties:
\newcounter{prop}
\begin{list}%
{\roman{prop})}{\usecounter{prop}\setlength{\leftmargin}{0.3in}}
\item\label{inc}
$\{M_k\}_{k\geq 20}$ is strictly increasing and
$M_{k+1}-M_k\to +\infty$ as $k\to\infty$.
\item\label{diff}
$\{M_k x_k\}_{k\geq 20}$ is strictly
increasing and $M_{k+1}x_{k+1}-M_kx_k\to +\infty$
as $k\to\infty$.
\item\label{cover}
For any integer $l\geq 1$, $s\in\R_+$, $\delta>0$, and $K>0$,
there exists an integer $k>K$ such that $j(k)=l$ and
$|s-(x_k+\frac{\alpha_l}{M_k})|<\frac{\delta}{M_k}$.
\end{list}

Returning to the notation of \ref{mi},
we set $r_k=e^{-M_kx_k}$.
Here the properties (i) and (ii) are necessary in order for a vector
from Step (\ref{vect}) of \ref{mi} to exist; the property (iii) is what will
translate into the condition ($2$) of the Lemma.

\subsection{Definition of $M_k$, $x_k$}\label{mainconst}
We will inductively construct  the sequence $\{M_k\}_{k\geq 20}$ 
simultaneously with
an auxilliary sequence of real numbers \newline
$\{y_k\}_{k\ge 20}$.
Later we will set $x_k=y_k-\frac{\alpha_{j(k)}}{M_k}$.

{\bf Step 1.} Set $M_{20}=1$, $y_{20}=1$.

Now we turn to the description of   the step number $q$ ($q>1$).
 In every step several consecutive terms $M_n$ and $y_n$ are defined;
denote the number of the last term defined in the first $(q-1)$ steps
by N. We shall always  assume that $j(N)=1$, that is,
a step can end only when the function $j$ returns to the value $1$.

{\bf Step $q$ .}
Set
\begin{equation}\label{jump}
M_{N+1}=q^2(M_N+1),\qquad y_{N+1}=\frac2q.
\end{equation}
Denote the difference $M_{N+1}-M_N$ by $d$.
Now let $N_1$ be the minimal integer such that $N_1>N$ and $j(N_1)=1$.
Then we set
\begin{equation}
\begin{align}
M_{k+1}&=M_k+d+q[\ln\ln k], &\qquad k=N+1,N+2,\ldots,N_1-1;\label{ms}\\
y_{k+1}&=y_k,&\qquad k=N+1,N+2,\ldots,N_1-2;\label{yk}\\
y_{N_1}&=y_{N_1-1}+\frac 1{qM_{N_1}}\label{ystep}.
\end{align}
\end{equation}
After that
 we repeat the procedure: find $N_2$ -- the minimal integer greater than
$N_1$ such that $j(N_2)=1$, and define $M_k$ by the formula
(\ref{ms}) for $k=N_1,N_1+1,\ldots,N_2-1$, while the $y$'s are kept equal
to $y_{N_1}$, until at $k=N_2$ the sequence $\{y_k\}$
again ``makes a step'' of
the length $\frac 1{qM_{N_2}}$.
In this fashion,
we proceed to construct the terms with the indices up to $N_3$,
$N_4$, and so on until $y_{N_m}$ becomes greater than $q$ for some $m$.
Then the Step number $q$ ends, so $M_{N_m}$ and $y_{N_m}$
are the last terms constructed on the Step $q$.
Notice that by construction $j(N_m)=1$ as required.
Below we will also need  the observation that in any case
$y_{N_m}<2q$, since every time $y$ increases,
it increases by a number less than $1$,
and we stop as soon as $y_{N_m}$ becomes greater than $q$.

This would have completed the construction, but we need to show that
the above procedure eventually terminates, that is,
 that $y_k$ eventually does become
greater than $q$. Here is a proof:

\subsection{Proof of  termination of the algorithm}\label{itends!}
The proof of 
 existence of an integer  $N_m$ such that $y_{N_m}>q$
is carried out by contradiction.
Suppose that $y_{N_i}\le q$ for $i=1,2,\dots$.
 By definition of the sequence
$y_k$, this means that the series
\begin{equation}\label{riad}
\sum_{{k>N, \,}\atop{j(k)=1}}\frac 1{qM_k}
\end{equation}
converges.
By the construction of the sequence $M_k$, we have
$M_k\le M_{N+1}+d(k-N-1)+q(k-N-1)\ln\ln(k-1)$
for $k>N$.
This implies that $M_k\le C k\ln\ln k$ for some constant $C$ that
depends on $q$ and $N$ but does not depend on $k$.

By definition of the function $j(n)$, the convergence of the series
(\ref{riad}) is equivalent to the convergence of the series
$\sum_{s\ge 20}\frac 1{M_{[s\ln\ln s]}}$ (recall that at the moment
we are assuming that the Step number $q$ never terminates, that is,
an infinite number of $M_k$'s
are defined by the formula (\ref{ms}), so this series is supposedly
 well defined).
Combining this statement with the estimate on $M_k$, we see that
our assumption ultimately implies that
the series
\begin{equation}\label{lnlnln}
\sum_{s\geq 20}\frac 1{s\cdot \ln\ln s\cdot\ln\ln(s\ln\ln s)}
\end{equation}
converges, which is not the case, and we arrive at a contradiction.
\qed

The construction of the sequences $\{M_k\}_{k\geq 20}$ and $\{y_k\}_{k\geq 20}$
is now completed.
Set
$$
x_k=y_k-\frac{\alpha_{j(k)}}{M_k}, \qquad k\geq 20.$$

\subsection{Verification of the properties}
We now turn to the verification of the required properties of these sequences.

The property (i) is immediate from the construction.
\subsubsection{Property (ii)}\label{prdva}
We need
an  estimate  from below for the difference $M_{k+1}x_{k+1}-M_k x_k$.
In view of the way in which the
sequences $\{M_k\}_{k\geq 20}$, $\{x_k\}_{k\geq 20}$
were constructed, there are two cases.

{\it Case $1$}.
The index $k$ is such that both $M_k$ and $M_{k+1}$ (and, accordingly, $x_k$
and $x_{k+1}$) are defined within the same step.

We will treat this case keeping the notations introduced in the description
of Step number $q$, $q>1$.
Since $y_{k+1}\geq y_k$ within Step $q$, we have
$$M_{k+1}x_{k+1}-
M_{k}x_{k}\ge
\\ y_{k}(d+q[\ln\ln k])-\alpha_{j(k+1)}+\alpha_{j(k)}.
$$
Recall that $j(k)<[\ln\ln k]+3$ by
(\ref{jin}) and $|\alpha_{l+1}-\alpha_l|\le 1$ by (\ref{alfa}).
These estimates together yield
$$
\alpha_{j(k)}-\alpha_{j(k+1)}=
\sum_{i=j(k)}^{j(k+1)-1}\alpha_{i}-\alpha_{i+1}\ge -(j(k+1)-j(k))>
-[\ln\ln(k+1)]-2.
$$
Combining all of the above with
the inequality $d\ge q^2$ which follows from (\ref{jump}), we get
$$M_{k+1}x_{k+1}-
M_{k}x_{k}>\frac 2 q (q^2+q[\ln\ln k])-[\ln\ln (k+1)] - 2> 2q-3.
$$

{\it Case $2$}. The index  $k$ happened to
be a ``boundary'' index, that is,  $M_k$ is defined on the Step q-1,
whereas $M_{k+1}$ is the first element of the $q$-th step.

Assume that $q>2$ (for $q=2$, we have $M_{21}x_{21}-M_{20}x_{20}>0$).
Again, we keep the notation of the previous section. Say,
$k=N$ -- the last
index of the Step q-1.
Then, by definiton (\ref{jump}),
$M_{k+1}=q^2 (M_k+1)$, and
 $x_{k+1}=\frac2q-\frac{\alpha_{j(k+1)}}{M_{k+1}}$.
Recall that
$j(N)=1$, $j(N+1)\le 2$ by definition of the function $j$, and
$\alpha_1=\alpha_2=0$, see (\ref{alfa}).
Besides, we know that $q-1<y_N<2(q-1)$ since $y_N$ was constructed on the 
Step q-1.
Finally,
\begin{eqnarray*}
M_{k+1}x_{k+1}-M_kx_k> \frac 2q M_{k+1}-2qM_k\\
\ge 2\left(\frac1q M_{k+1}-qM_k\right)=
2\left(\frac{q^2M_k+q^2}q-qM_k\right)=2q.
\end{eqnarray*}

We see that $M_{k+1}x_{k+1}-
M_{k}x_{k}$ is always positive and increases to $\infty$.\qed

\subsubsection{Property (iii)}
Let us fix the numbers $s\in\R_+$, $\delta>0$, and $l, K\in\N$.
Pick a number $Q>K$ such that $\frac1Q<\delta$ and $\frac2Q< s <Q$.
Let $\{n_m\}_{m\ge 1} $ be the set of all positive integers $n$, such that
$j(n)=l$.
We will think of $\{n_m\}$ as an increasing sequence of positive integers.
Observe that all the numbers $M_{n_m}$ with sufficiently large $m$ were
defined on the steps with numbers $q>Q$.
Let us consider all the  points
$y_{n_m}=x_{n_m}+\frac{\alpha_l}{M_{n_m}}$
which were constructed on the $q$-th step.
Their
$\frac1{qM_{n_m}}$-neighbourhoods
cover the interval
$]\frac2q,\,q[$ by construction, and since the
$\frac{\delta}{M_{n_m}}$-neighbourhoods are even larger, $s$ falls in one
of them.\qed

\section{Common hypercyclic vector}
\subsection{Construction of the vector}\label{real}
As announced in the previous section, we set $r_k=e^{-M_kx_k}$.
By the property (ii), the sequence $r_k$ is decreasing and it
tends to $0$ as $k\to\infty$.

We turn to the construction of a vector $f\in \ldva$ such that
$\|B^{M_k}f\|=r_k$ and $B^{M_k}f$ belongs to the cone $C_{j(k)}$
for all sufficiently large $k\in\N$.

Let $\{v_l\}_{l\geq 1}$ be the dense subset of $\ldva$ from
\ref{vectorsandcones}.
For an integer $k\geq 20$, let $d_k=\sqrt{r_k^2-r_{k+1}^2}$, and
let $w_k\in\ldva$ be the vector whose first $M_{k+1}-M_k$ coordinates
coincide with those of $v_{j(k)}$ and the coordinates starting from
$M_{k+1}-M_k+1$
equal  $0$.
Then set
\begin{equation}
f=\sum_{k\geq 20}\frac {d_k}{\|w_k\|}S^{M_k}w_k,
\end{equation}
where
$S$  is the forward shift operator,
$S(a_0,a_1,a_2\ldots)=(0,a_0,a_1,\ldots)$.

It is clear that
$\|B^{M_k} f\|=\left(\sum_{i=k}^{\infty}d_i^2\right)^{1/2}=r_k$.

We now turn to show that for any $l\ge 1$, the vector $B^{M_k}f$ falls into
the
cone $\co_l$ for all sufficiently large $k$ such that $j(k)=l$. First, observe
that
$
\|B^{M_k}f-\frac {d_k}{\|w_k\|}w_k\|=r_{k+1}
$
by the definition of $f$.
Second, the property (ii) of $M_k$, $x_k$ implies that
$\frac{d_k^2}{r_{k+1}^2}=\frac{r_k^2-r_{k+1}^2}{r_{k+1}^2}\to+\infty$.
Combining these two facts, we get
\begin{equation}\label{last}
 \left\|\frac{B^{M_k}f}{d_k}-\frac {w_k}{\|w_k\|}\right\|=
\frac{r_{k+1}}{d_k}\to 0 \,\,\,{\text {\ as \ }}\, k\to\infty.
\end{equation}

The vectors $v_k$ were assumed to be finite.
 Therefore, for a fixed $l$, for all sufficiently large $k$
such that $j(k)=l$, the vectors $v_l$ and $w_k$ coincide, since
$M_{k+1}-M_k\to\infty$ as $k\to\infty$.
Thus $(\ref{last})$ yields
$\|\frac{\|v_l\|}{d_k}B^{M_k}f-v_l\|<\epsilon_l$, and hence
$B^{M_k}f$ belongs to the cone $\co_{j(k)}$
for all sufficiently
large $k$ such that $j(k)=l$.

\subsection{Proof of hypercyclicity for {\bf $\lambda>1$}}
To show that $f$ is a common hypercyclic vector for all operators
$\lambda B$ with $\lambda>1$, it suffices  to verify that for any fixed
 integer $l\ge 1$  there exists $k\geq 20$ such that
$j(k)=l$ and
\begin{equation}\label{hc}
\|(\lambda B)^{M_k} f-v_l\|<3\epsilon_l.
\end{equation}

Fix $l\ge 1$, and
let $K$ be  a positive  integer such that for all $n>K$ the vector $B^{M_n}f$
falls into the cone $\co_l$ whenever $j(n)=l$.
The existence of $K$ was proved in
the previous subsection.
Now the statement immediately follows from Lemma.
Indeed,
 the property (iii) implies that the condition
($2$) of the Lemma  is satisfied if we set $a=e^{\alpha_l}
=\|v_l\|$ and let
$\{n_k\}_{k\geq 1}$
 be the sequence of all numbers $M_n$ with $n>K$ and $j(n)=l$.
This is precisely one of the subsequences discussed in Step 4 of the 
algorithm described in \ref{mi}.
Then by the Lemma
for any $\lambda>1$ there exists $k>K$ with $j(k)=l$ such that
$|\lambda^{M_k}\|B^{M_k}f\|-\|v_l\| |<\epsilon_l$ (recall 
that $r_k = \|B^{M_k}f\|$).
Since the vector $(\lt)^{M_k}f$
belongs to the cone generated by the ball of radius
$\epsilon_l$ centered at $v_l$, and its norm differs from the norm of $v_l$
by less than $\epsilon_l$, the condition (\ref{hc}) is fulfilled.

\qed

\subsection{Complex case}\label{complex}
Let us now consider the complex multiples $zB, \, |z|>1$,  of the operator $B$.
We write $z$ in the trigonometric
form $z=\lambda e^{2\pi i\theta}$ with $\lambda >1$ and $\theta\in [0,1[$.

We preserve the same notation as before, except for the function $j(k)$
which should be defined here by means of the sequence $n \ln\ln\ln n$ (and not
$n\ln\ln n$ as in Section 1.3.2).
We fix a cone $\co_l$, exactly as before.
For real $\lambda$'s it was
necessary to establish that $\frac{\delta}{M_k}$-neighbourhoods
of the points $\frac{\ln\|v_l\|-\ln\|B^{M_k} f\|}{M_k}$ (with $M_k$ such that
$B^{M_k}f$ belongs to $\co_l$) cover
$\R_+$ for any $\delta>0$. We achieved this by
constructing a ``very dense'' sequence $\{x_k\}  $ in $\R_+$ (see
Section \ref{mainconst}) and then
 producing a vector $f$ such that
$x_k=-\ln\|B^{M_k} f\|$.

Now, in the same fashion as was previously done for the sequence
 $\{x_k\}$, 
we construct the sequence
$\{(x_k,\theta_k)\}$ of elements of the set $\R_+\times[0,1[$.
More precisely, let $y_k$ be as in Section \ref{mainconst};
we will be constructing $(y_k,\theta_k)$ simultaneously with $M_k$.
Let us look back at ``Step q''.
We keep defining $y_k$'s and $M_k$'s by the same formulas.
The only modification is that now before $y_k$ is allowed
to make a step forward (see the formula (\ref{ystep})), the
$\theta$-component will have to increase from $0$ to $1$ with the
decreasing  steps $\frac 1{qM_{N_i}}$, also staying constant
at its every value
with only $M$'s changing until the function $j(n)$ completes a
 cycle and returns to the value $1$.
The same argument as in Section \ref{itends!} shows that $y_k$ 
is allowed to move forward each
time once $\theta$ has completed the cycle.
The proof that each step terminates essentially remains the same.

Given the sequences $\{M_k\}$ and $\{(x_k,\theta_k)\}$,
we define the vector $f$ by
\begin{equation}
f=\sum_{k
}\frac {d_k}{\|w_k\|}e^{2\pi i\theta_k}
(S^{M_k}w_k)
\end{equation}
in the notation of \ref{real}.
The proof of its hypercyclicity works in the same way as
in the previous section, with the use of the
fact that now the square $\frac{\delta}{M_k}$-neighbourhoods
of the points $(y_k,\theta_k)$ cover the half strip $\R_+\times [0,1[$.

\vskip 0.3cm
In conclusion, we remark that 
the theorem remains true (with essentially the same proof) for the backward
shift acting on the
spaces $\ell^p, \, 1\le p<\infty$, and on the space $c_0$ of sequences 
that tend to $0$.

%
%
%
%

\end{document}